\documentclass[12pt]{amsart}
\usepackage{latexsym}
\usepackage{amssymb,times,fullpage}
\newtheorem{thm}{Theorem}

\newtheorem{lem}[thm]{Lemma}

\theoremstyle{remark}
\newtheorem{rem}[thm]{Remark}

\theoremstyle{definition}

\newcommand{\R}{\mathbb{ R}}

\newcommand{\Z}{\mathbb{ Z}}

\title[Failure of separation in mapping class groups]{Failure of 
separation by quasi-homomorphisms\\ in mapping class groups}
\author{H.~Endo}
\address{Department of Mathematics, Graduate School of Science, Osaka University, 
Toyonaka, Osaka 560-0043, Japan}
\email{endo@math.wani.osaka-u.ac.jp}
\author{D.~Kotschick}
\address{Mathematisches Institut, Ludwig-Maximilians-Universit\"at M\"unchen,
Theresienstr.~39, 80333 M\"unchen, Germany}
\email{dieter@member.ams.org}
\thanks{The second author would like to thank L.~Polterovich for a conversation raising the question whether 
a separation theorem for mapping class groups of higher genus surfaces holds, and K.~Fujiwara and J.~McCarthy 
for useful comments. Support from the {\it Deutsche Forschungsgemeinschaft} and from JSPS Grant 18540083
is gratefully acknowledged}
\date{June 2, 2006; MSC 2000: primary 20F65, secondary 20F12, 20F69, 57M07}

\begin{document}

\begin{abstract}
    We show that mapping class groups of surfaces of genus at least 
    two contain elements of infinite order that are not conjugate 
    to their inverses, but whose powers have bounded torsion lengths. 
    In particular every homogeneous quasi-homomorphism vanishes on 
    such an element, showing that elements of infinite order not conjugate 
    to their inverses cannot be separated by quasi-homomorphisms.
\end{abstract}

\maketitle

This note was motivated by the work of Polterovich and 
Rudnick~\cite{PR,PR2}, who used the geometry of the hyperbolic plane 
to show that on $SL_{2}(\Z)$ quasi-homomorphisms exist in abundance 
and have interesting properties. One of their results is the following:
\begin{thm}[\cite{PR}]\label{t:PR1}
    Let $g\in SL_{2}(\Z)$ be an element of infinite order not 
    conjugate to its inverse. Then there exists a homogeneous 
    quasi-homomorphism $\varphi\colon SL_{2}(\Z)\rightarrow\R$ with 
    $\varphi(g)\neq 0$.
    \end{thm}
Since homogeneous quasi-homomorphisms satisfy 
$\varphi(g^{-1})=-\varphi(g)$ and are constant on conjugacy classes, 
the assumption of the theorem is clearly necessary\footnote{The reader
can consult~\cite{Ko} for background on quasi-homomorphisms.}.
    
    In view of Brooks's classical constructions of quasi-homomorphisms, 
this result is not suprising. In fact, Polterovich and 
Rudnick~\cite{PR} pointed out that such a result essentially holds in 
much greater generality, because for all non-elementary Gromov hyperbolic 
groups it can be deduced from the work of Epstein and 
Fujiwara~\cite{EF}.

In~\cite{PR2}, Polterovich and Rudnick generalized Theorem~\ref{t:PR1} 
to the following ``separation theorem'':
\begin{thm}[\cite{PR2}]\label{t:PR2}
    Let $g\in SL_{2}(\Z)$ be a primitive element of infinite order not 
    conjugate to its inverse, and $g_{1},\ldots,g_{n}\in SL_{2}(\Z)$ 
    any finite number of elements not conjugate to any power of $g$.
    Then there exists a homogeneous quasi-homomorphism 
    $\varphi\colon SL_{2}(\Z)\rightarrow\R$ with $\varphi(g)\neq 0$ 
    and $\varphi(g_{1})=\ldots=\varphi(g_{n})=0$.
    \end{thm}

Thinking of $SL_{2}(\Z)$ as the mapping class group of the two-torus, one 
naturally wonders whether Theorems~\ref{t:PR1} and~\ref{t:PR2} can 
be generalized to the mapping class groups of higher-genus surfaces. 
Several years ago we proved that there are non-trivial homogeneous
quasi-homomorphisms on mapping class groups~\cite{EK}. Bestvina and 
Fujiwara~\cite{BF} then showed that the space of such quasi-homomorphisms 
is infinite-dimensional, and Polterovich asked whether it might be possible to 
prove a separation theorem in the spirit of Theorem~\ref{t:PR2} for mapping 
class groups. On the one hand, mapping class groups are perfect if the 
genus of the underlying surface is at least three~\cite{Powell}, so
they certainly have no homomorphisms to Abelian groups. On the other 
hand, they are residually finite~\cite{Gr}, thus their elements can be 
separated by homomorphisms to finite groups, and by linear 
representations. 

It is our purpose here to show that elements in mapping class groups cannot 
be separated by quasi-homomorphisms, by showing that the analogue of 
Theorem~\ref{t:PR1} and, {\it a fortiori}, the analogue of Theorem~\ref{t:PR2} 
fail for mapping class groups of surfaces of genus $\geq 2$. We shall prove 
the following:
\begin{thm}\label{t:main}
    For every closed oriented surface of genus at least $2$ there exist 
    primitive elements $g$ of infinite order in its mapping class group 
    of orientation-preserving diffeomorphisms such that $g^k$ is not 
    conjugate to $g^{-k}$ for all $k\neq 0$, but all powers of $g$ are 
    products of some fixed number of torsion elements.
    \end{thm}
It follows from the boundedness of the torsion lengths $t(g^{n})$ that 
the stable torsion length 
$$
\vert\vert 
g\vert\vert_{T}=\lim_{n\rightarrow\infty}\frac{t(g^{n})}{n} 
$$
vanishes, which, by the results of~\cite{Ko}, implies the vanishing of 
the stable commutator length. In fact, we will check explicitly that the 
powers of $g$ have bounded commutator lengths. The vanishing of the stable 
torsion length of $g$ also implies that every homogeneous quasi-homomorphism 
must vanish on $g$, by the estimate
$$   
\vert\vert 
g\vert\vert_{T}\geq\frac{\vert\varphi(g)\vert}{D(\varphi)} \ ,
$$
where $D(\varphi)$ denotes the defect of $\varphi$, compare~\cite{Ko}. 

To put Theorem~\ref{t:main} into perspective, recall that, on the one hand, mapping 
class groups are not hyperbolic, for example because they contain 
Abelian subgroups of large ranks generated by Dehn twists along 
disjoint curves, cf.~\cite{BLM}. Reducible elements that are products 
of commuting Dehn twists, or have powers which are such products, 
are crucial to our construction, which is a further generalization of 
examples in~\cite{Ko}, with roots in the work of McCarthy and 
Papadopoulos~\cite{MP}. On the other hand, the construction of 
quasi-homomorphisms due to Bestvina and Fujiwara~\cite{BF} does not 
require the group to be hyperbolic. All that is needed is a suitably 
weakly proper action on a $\delta$-hyperbolic space, and this 
weaker property holds for the actions of mapping class groups on 
complexes of curves, because of the work of Masur and Minsky~\cite{MM}. 
Polterovich and Rudnick~\cite{PR2} suggest that the techniques of~\cite{BF} 
might lead to a separation theorem for hyperbolic groups. Although we 
show that the separation theorem cannot hold for weakly hyperbolic 
groups like the mapping class groups, it is still possible that a 
separation theorem does hold if one restricts to pseudo-Anosov elements 
of mapping class groups, because only these act hyperbolically on the 
curve complex, see~\cite{MM,BF}. Results in this direction are contained 
in forthcoming work of Calegari and Fujiwara~\cite{CF}.
Their results are in contrast to Theorems~\ref{t:PR1} 
and~\ref{t:PR2}, where $g$ is not assumed to be hyperbolic.

\begin{proof}[The proof of Theorem~\ref{t:main}]
We denote by $\Sigma=\Sigma_{h}$ a fixed closed oriented surface of genus 
$h\geq 2$. Let $a$, $b$ and $c$ be disjoint non-separating simple closed 
curves on $\Sigma$ in distinct isotopy classes, and let 
$f=t_{a}^{-1}t_{b}^{-1}t_{c}^{2}$, where $t_{\alpha}$ denotes the 
right-handed Dehn twist along $\alpha$. Then 
$f^{k}=t_{a}^{-k}t_{b}^{-k}t_{c}^{2k}$, for all integers $k$.
In particular, $f$ is of infinite order in the mapping class group 
of $\Sigma$.

Now we use the following obvious lemma:
\begin{lem}\label{l:l}
    Let $a$ and $b$ be non-separating disjoint simple closed curves on $\Sigma$. 
    Then there exists an involution $\varphi\colon\Sigma\rightarrow\Sigma$ 
    interchanging $a$ and $b$.
    \end{lem}
This allows us to find involutions $\varphi$ interchanging 
$a$ and $c$, and $\psi$ interchanging $b$ and $c$. Then 
$$
f^{k}=t_{c}^{k}t_{a}^{-k}t_{c}^{k}t_{b}^{-k}=t_{c}^{k}\varphi 
(t_{c}^{k})^{-1}\cdot\varphi\cdot t_{c}^{k}\psi 
(t_{c}^{k})^{-1}\cdot\psi = [t_{c}^{k},\varphi]\cdot [t_{c}^{k},\psi] \ .
$$
Thus $f^{k}$ is both a product of four involutions and a product of two commutators.
Therefore both the torsion length and  the commutator length of $f^{k}$ are
bounded independently of $k$. In particular, all homogeneous quasi-homomorphisms vanish 
on $f$.

Suppose now that $f^{k}$ and $f^{l}$ are conjugate for some $k$ and $l$. 
The system of curves consisting of $a$, $b$ 
and $c$ is a reducing system for all powers of $f$. In fact, in 
the terminology of Matsumoto and Montesinos~\cite{MatsMon}, these 
three curves form a precise reducing system. As precise reducing 
systems are unique up to isotopy, a conjugacy between $f^{k}$ and 
$f^{l}$ must permute the isotopy classes of these three curves in 
such a way that the twisting numbers along the curves are matched up. 
But the triples of twisting or skrew numbers are $(-k,-k,2k)$ and 
$(-l,-l,2l)$, so that we conclude $k=l$. In particular, $f$ cannot be 
conjugate to its inverse.

The element $f$ considered above is not always primitive, as it may be 
expressed as the square of $\alpha t_{a}^{-1}t_{c}$, for any involution 
$\alpha$ interchanging $a$ and $b$ and fixing $c$ (not necessarily pointwise). 
However, as the mapping class group of $\Sigma$ is residually 
finite~\cite{Gr}, there exists a primitive element $g$ such that $f$ is a 
power of $g$. Then all powers of $g$ also have bounded torsion and commutator 
lengths, and every homogeneous quasi-homomorphism vanishes on $g$. If two 
distinct powers of $g$ were conjugate to each other, then the same would be 
true for $f$.

This completes the proof of Theorem~\ref{t:main}.
\end{proof}

\begin{rem}
The ideas on reducing systems from~\cite{MatsMon} that we have used to show
that distinct powers of $f$ are not conjugate to each other are quite classical,
going back at least to Nielsen. Another modern approach to this topic is due to Birman, Lubotzky
and McCarthy~\cite{BLM}, and we could just as well base our argument on this 
reference.
\end{rem}

\begin{rem}\label{r}
    There are of course many other reducible mapping classes to which similar arguments apply. 
    For example, assuming $h\geq 3$, we can find disjoint non-separating non-isotopic simple 
    closed curves $a_{1}$, $a_{2}$, $b_{1}$ and $b_{2}$ with the following 
    property: 
    \begin{itemize}
	\item[($\star$)] \ \ \ \ $\Sigma\setminus (a_{1}\cup a_{2})$ is connected, but 
    $\Sigma\setminus (b_{1}\cup b_{2})$ is not.
    \end{itemize}
    By Lemma~\ref{l:l}  there exist involutions $\varphi_{i}$ interchanging $a_{i}$ 
    with $b_{i}$.
    The element $g=t_{a_{1}}t_{a_{2}}t_{b_{1}}^{-1}t_{b_{2}}^{-1}$ has the property that 
    $$
    g^{k}=t_{a_{1}}^{k}t_{b_{1}}^{-k}t_{a_{2}}^{k}t_{b_{2}}^{-k} = 
    t_{a_{1}}^{k}\varphi_1 t_{a_{1}}^{-k}\cdot \varphi_1 \cdot t_{a_{2}}^{k}\varphi_2t_{a_{2}}^{-k}\cdot\varphi_2
    = [t_{a_{1}}^{k},\varphi_1]\cdot [t_{a_{2}}^{k},\varphi_2] \ ,
    $$
    which is a product of 4 involutions, and also a product of 2 commutators. 
    
    In this case the set of curves $a_{1}$, $a_{2}$, $b_{1}$ and 
    $b_{2}$ forms a precise reducing system for every $g^{k}$, and the 
    twisting numbers are $(k,k,-k,-k)$. Now if $g^k$ were conjugate to 
    $g^{-k}$, then a conjugacy would have to interchange the isotopy 
    classes of $a_{1}\cup a_{2}$ and of $b_{1}\cup b_{2}$, but this is impossible 
    by ($\star$).
    \end{rem}

\begin{rem}
    If an element $g$ in a mapping class group can be expressed as a 
    product of positive powers of Dehn twists along disjoint curves, 
    then by~\cite{EK,Ko} there is a homogeneous quasi-homomorphism $\varphi$ 
    with $\varphi(g)\neq 0$. On the other hand, it follows from Lemma~\ref{l:l} that 
    $t_{a}t_{b}^{-1}$ is conjugate to its inverse for any two disjoint non-separating 
    simple closed curves. Thus the examples we have given in the proof of 
    Theorem~\ref{t:main} and in Remark~\ref{r} are the simplest reducible mapping 
    classes which exhibit the failure of Theorems~\ref{t:PR1} and~\ref{t:PR2} for 
    mapping class groups. 
\end{rem}

\begin{rem}
Our examples show that the linear growth of the commutator length for products
of positive powers of Dehn twists along disjoint curves depends on restricting the 
chirality of the Dehn twists, even after excluding elements conjugate to their inverses. 
For linear growth of the word length no such restriction of the chirality is necessary~\cite{FLM}.
\end{rem}

\bigskip

\bibliographystyle{amsplain}

\bigskip

\end{document}